\documentclass[12pt,leqno]{article}

\usepackage{amsmath,amsthm,amsfonts,latexsym,amscd,amssymb}
\numberwithin{equation}{section}
\input xypic
\def\qed{{\hbadness=10000\hfill\ \vbox{\hrule height.09ex
   \hbox{\vrule width.09ex height1.55ex depth.2ex \kern1.8ex
   \vrule width.09ex height1.55ex depth.2ex}\hrule height.09ex}\break
   \bigskip}}
\newtheorem{theorem}{Theorem}[section]
\newtheorem{lemma}{Lemma}[section]

\theoremstyle{definition}

\theoremstyle{remark}

\begin{document}

\linespread{1}\title{\textbf{Hypersurface of a Finsler space subjected to a Kropina change with an \textsl{h}-vector}}

\author{M.$\,$K. \textsc{Gupta}\thanks{Corresponding author. E-mail: mkgiaps@gmail.com}~\\
\normalsize{Department of Pure $\&$ Applied Mathematics}\\
\normalsize{Guru Ghasidas Vishwavidyalaya}\\
\normalsize{Bilaspur (C.G.), India}\\ \\
P.$\,$N.$\,\,$\textsc{Pandey}\\
\normalsize{Department of Mathematics}\\
\normalsize{University of Allahabad}\\
\normalsize{Allahabad, India}}
\date{}
\maketitle

\linespread{1.3}\begin{abstract} The concept of \textsl{h}-vector was introduced by H. Izumi in 1980. Recently we have obtained the Cartan connection for the Finsler space whose metric is given by Kropina change with an \textsl{h}-vector. In 1985, M. Matsumoto studied the theory of Finsler hypersurface. In this paper, we derive certain geometrical properties of a Finslerian hypersurface subjected to a Kropina change with an \textsl{h}-vector.
\end{abstract}

\section{Introduction}
~~~~In 1984, C. Shibata \cite{Sh84}\linespread{1.3}\footnote{Numbers in square brackets refer to the references at the end of the paper.} dealt with a change of Finsler metric which is called a $\beta$-change of metric. A remarkable class of $\beta$-change is Kropina change $^{*}L(x,y)=\frac{L^2(x,y)}{b_i(x)\,y^i}$. If $L(x,y)$ is a metric function of a Riemannian space then $^{*}L(x,y)$ reduces to the metric function of a Kropina space. Such a Finsler metric was introduced by V.K. Kropina \cite{Kr61}.

H. Izumi \cite{Iz80} while studying the conformal transformation of Finsler spaces, introduced the \textsl{h}-vector $b_{i}$, which is \textsl{v}-covariant constant with respect to the Cartan connection and satisfies $L\,C^{h}_{ij}\,b_{h}=\rho\, h_{ij}\,,\,\rho\neq 0$. Thus if $b_{i}$ is an \textsl{h}-vector then
\begin{equation}(i)\,\,b_{i}|_{k}=0\,,\qquad\qquad\qquad(ii)\,\, L\,C^{h}_{ij}\,b_{h}=\rho\,h_{ij}\,.\end{equation}
This gives
\begin{equation}L\,\dot{\partial_{j}}b_{i}=\rho \,h_{ij}\,. \end{equation}
Since $\rho\neq 0$ and $h_{ij}\neq0$ for $h_{ij}=0$ implies $n=1$ which is not true, the \textsl{h}-vector $b_{i}$ depends not only on positional coordinates but also on directional arguments. Izumi \cite{Iz80} showed that $\rho$ is independent of directional arguments. The present authors \cite{GuPa13} obtain the relation between the Cartan connections of $F^{n}=(M^{n},L)$ and $^{*}F^{n}=(M^{n},\,^{*}L)$ where $^{*}L(x,y)$ is obtained by the transformation
\begin{equation}^{*}L(x,y)=\frac{L^2(x,y)}{b_i(x,y)\,y^i}\end{equation}
and $b_{i}(x,y)$ is an \textsl{h}-vector in $(M^{n},L)$.

M. Matsumoto \cite{Ma85} presented a systematic theory of Finslerian hypersurface. The present authors \cite{Ma08} obtained certain results for the Finslerian hypersurface.

In this paper, certain geometrical properties of a Finslerian hypersurface subjected to a Kropina change with an \textsl{h}-vector, were disscussed.

The terminologies and notations are referred to Matsumoto \cite{Ma86}.

\section{Preliminaries}
~~~~~Let $M^{n}$ be an $n$-dimensional smooth manifold and $F^{n}=(M^{n},L)$ be an $n$-dimensional Finsler space equipped with a metric function $L(x,y)$ on $M^{n}$. The normalized supporting element, the metric tensor, the angular metric tensor and Cartan tensor are defined by $l_{i}=\dot{\partial_{i}}L\,,\,\, g_{ij}=\frac{1}{2}\,\dot{\partial_{i}}\,\dot{\partial_{j}}L^{2}\,,\,\, h_{ij}=L\,\dot{\partial_{i}}\,\dot{\partial_{j}}L\,$ and  $\,C_{ijk}=\frac{1}{2}\,\dot{\partial_{k}}\,g_{ij}\,$ respectively. Throughout this paper, we use the symbols $\dot{\partial_{i}}$ and $\partial_{i}$ for $\partial/\partial y^{i}$ and $\partial/\partial x^{i}$ respectively. The Cartan connection in $F^{n}$ is given as $C\Gamma=(F^{\,i}_{jk},G^{\,i}_{j},C^{\,i}_{jk})$. The \textsl{h}- and \textsl{v}-covariant derivatives of a covariant vector $X_{i}(x,y)$ with respect to the Cartan connection are given by
\begin{equation}X_{i|j}=\partial_{j}X_{i}-(\dot{\partial}_{h}X_{i})G^{h}_{j}-F^{\,r}_{ij}X_{r}\,,\end{equation}\\[-11mm]
and\\[-11mm]
\begin{equation} X_{i}|_{j}=\dot{\partial}_{j}X_{i}-C^{\,r}_{ij}X_{r}\,.\end{equation}

A hypersurface $M^{n-1}$ of the underlying manifold $M^{n}$ may be represented parametrically by the equations $x^{i}=x^{i}(u^{\alpha})$, where $u^{\alpha}$ are the Gaussian coordinates on $M^{n-1}$ (Latin indices run from 1 to $n$, while Greek indices take values from 1 to $n$-1). We assume that the matrix of projection factors $B^{i}_{\alpha}=\partial x^{i}/\partial u^{\alpha}$ is of rank $n$-1. If the supporting element $y^{i}$ at a point $u=(u^{\alpha})$ of $M^{n-1}$ is assumed to be tangent to $M^{n-1}$, we may then write $y^{i}=B^{i}_{\alpha}(u)\,v^{\alpha}$ so that $v=(v^{\alpha})$ is thought of as the supporting element of $M^{n-1}$ at the point $u^{\alpha}$. Since the function $\underline{L}\,(u,v)=L\Big(x(u),y(u,v)\Big)$ gives rise to a Finsler metric on $M^{n-1}$, we get an ($n$-1)-dimensional Finsler space $F^{n-1}=(M^{n-1},\underline{L}\,(u,v))$.

At each point $u^{\alpha}$ of $F^{n-1}$, a unit normal vector $N^{i}(u,v)$ is defined by
\begin{equation} g_{ij}\,B^{i}_{\alpha}\,N^{j}=0\,,\qquad g_{ij}\,N^{i}\,N^{j}=1\,.\end{equation}
The inverse projection factors $B^{\alpha}_{i}(u,v)$ of $B^{i}_{\alpha}$ are defined as
\begin{equation} B^{\alpha}_{i}=g^{\alpha\beta}\,g_{ij}\,B^{j}_{\beta}\,,\end{equation}
where $g^{\alpha\beta}$ is the inverse of the metric tensor $g_{\alpha\beta}$ of $F^{n-1}$. \\From (2.3) and (2.4), it follows that
\begin{equation}B^{i}_{\alpha}\,B^{\beta}_{i}=\delta^{\beta}_{\alpha}\,,
\quad B^{i}_{\alpha}\,N_{i}=0\,,\quad N^{i}\,B^{\alpha}_{i}=0\,,\quad N^{i}\,N_{i}=1\,,\end{equation}
and further
\begin{equation}B^{i}_{\alpha}\,B^{\alpha}_{j}+N^{i}\,N_{j}=\delta^{i}_{j}\,.\end{equation}
For the induced Cartan connection $IC\Gamma=(F^{\alpha}_{\beta\gamma}\,,G^{\alpha}_{\beta}\,,C^{\alpha}_{\beta\gamma})$ on $F^{n-1}$, the second fundamental \textsl{h}-tensor $H_{\alpha\beta}$ and the normal curvature vector $H_{\alpha}$ are given by
\begin{equation}H_{\alpha\beta}=N_{i}\,(B^{i}_{\alpha\beta}+F^{i}_{jk}\,B^{j}_{\alpha}\,B^{k}_{\beta})+M_{\alpha}\,H_{\beta}\,,\end{equation}\\[-11mm]
and\\[-11mm]
\begin{equation}H_{\alpha}=N_{i}\,(B^{i}_{0\alpha}+G^{i}_{j}\,B^{j}_{\alpha})\,,\end{equation}
where $M_{\alpha}=C_{ijk}\,B^{i}_{\alpha}\,N^{j}\,N^{k}$, $B^{i}_{\alpha\beta}=\partial^{2}x^{i}/\partial u^{\alpha}\,\partial u^{\beta}$ and $B^{i}_{0\alpha}=B^{i}_{\beta\alpha}\,v^{\beta}$.\\
The equations (2.7) and (2.8) yield
\begin{equation}H_{0\alpha}=H_{\beta\alpha}\,v^{\beta}=H_{\alpha}\,,\quad\quad H_{\alpha 0}=H_{\alpha\beta}\,v^{\beta}=H_{\alpha}+M_{\alpha}\,H_{0}\,. \end{equation}
The second fundamental \textsl{v}-tensor $M_{\alpha\beta}$ is defined as:
\begin{equation} M_{\alpha\beta}=C_{ijk}\,B^{i}_{\alpha}\,B^{j}_{\beta}\,N^{k}.\end{equation}
The relative \textsl{h}- and \textsl{v}-covariant derivatives of $B^{i}_{\alpha}$ and $N^{i}_{}$ are given by
\begin{equation}\begin{split} B^{i}_{\alpha|\beta}=H_{\alpha\beta}\,N^{i}_{}\,,\quad\quad B^{i}_{\alpha}|_{\beta}=M_{\alpha\beta}\,N^{i}_{}\,,~~~~~~~\\
N^{i}_{\,|\beta}=-H_{\alpha\beta}\,B^{\alpha}_{j}\,g^{ij}\,,\quad\quad N^{i}_{}|_{\beta}=-M_{\alpha\beta}\,B^{\alpha}_{j}\,g^{ij}\,.\end{split}\end{equation}
Let $X_{i}(x,y)$ be a vector field of $F^{n}$. The relative \textsl{h}- and \textsl{v}-covariant derivatives of $X_{i}$ are given by
\begin{equation} X_{i|\beta}=X_{i|j}\,B^{j}_{\beta}+X_{i}|_{j}\,N^{j}\,H_{\beta}\,,\quad\quad X_{i}|_{\beta}=X_{i}|_{j}\,B^{j}_{\beta}\,.\end{equation}

Matsumoto \cite{Ma85} defined different kinds of hyperplanes and obtained their characteristic conditions, which are given in the following lemmas:
\begin{lemma}\label{lem1}
A hypersurface $F^{n-1}$ is a hyperplane of the first kind if and only if $H_{\alpha}=0$ or equivalently $H_{0}=0$. \end{lemma}
\begin{lemma}\label{lem2}
A hypersurface $F^{n-1}$ is a hyperplane of the second kind if and only if $H_{\alpha\beta}=0$. \end{lemma}
\begin{lemma}\label{lem3}
A hypersurface $F^{n-1}$ is a hyperplane of the third kind if and only if $H_{\alpha\beta}=0=M_{\alpha\beta}$. \end{lemma}

\section{The Finsler space $^{*}F^{n}=(M^{n},\,^{*}L)$}
~~~~If we denote $b^{}_{i}\,y^i$ by $\beta$ then indicatory property of $h_{ij}$ yield $\dot{\partial_{i}}\beta=b_i$. Throughout this paper, the geometric objects associated with $^{*}F^{n}$ will be asterisked. We shall use the notation $L_{ij}=\dot{\partial_{i}}\,\dot{\partial_{j}}L\,,\,\,L_{ijk}=\dot{\partial_{k}}\,L_{ij}\,,\ldots$ etc.
From (1.3), we get
\begin{equation}^{*}L_{ij}=(2\tau-\rho\tau^2)\,L_{ij}+\frac{2\,\tau^2}{\beta}\,m_i\,m_j\,,\end{equation}
\begin{equation}\begin{split}^{*}L_{ijk}=&(2\tau-\rho\tau^2)\,L_{ijk}+\frac{2\,\tau}{\beta}\,(\rho\tau-1)\,(m_iL_{jk}+m_jL_{ik}+m_kL_{ij})\\
&-\frac{2\tau^2}{L\beta}(m_im_jl_k+m_jm_kl_i+m_km_il_j)-\frac{6\tau^2}{\beta^2}m_im_jm_k\,,\end{split}\end{equation}
where $\tau=L/\beta\,,\,\,m_{i}=b_{i}-\frac{1}{\tau}l_{i}$. The normalized supporting element, the metric tensor of $^{*}F^{n}$ are obtained as \cite{GuPa13}
\begin{equation}^{*}l_{i}=2\,\tau l_{i}-\tau^2\,b_{i}\,,\end{equation}
\begin{equation}^{*}g^{}_{ij}=(2\tau^2-\rho\tau^3)\,g^{}_{ij}+3\tau^4\,b_{i}b_{j}-4\tau^3(l_{i}b_{j}+b_{i}l_{j})+(4\tau^2+\rho\tau^3)l_{i}l_{j}\,,\end{equation}
Differentiating the angular metric tensor $h_{ij}$ with respect to $y^k$, we get
\[\dot{\partial_{h}}h^{}_{ij}=2C^{}_{ijk}-\frac{1}{L}\,(l_i\,h_{jk}+l_j\,h_{ik}),\]
which gives
\begin{equation}L_{ijk}=\frac{2}{L}C^{}_{ijk}-\frac{1}{L^2}\,(h_{ij}l_k+h_{jk}l_i+h_{ki}l_j).\end{equation}
Taking (3.5) into account, (3.2) can be rewritten as
\begin{equation}^{*}C_{ijk}=(2\tau^2-\rho\tau^3)\,C_{ijk}-\frac{\tau^2}{2\,\beta}(4-3\rho\tau)(h_{ij}m_{k}+h_{jk}m_{i}+h_{ki}m_{j})-\frac{6\tau^2}{\beta}m_im_jm_k\,.\end{equation}
The inverse metric tensor of $^{*}F^{n}$ is derived as follows \cite{GuPa13}:
\begin{equation}\begin{split}^{*}g^{ij}=(2\tau^2-\rho\tau^3)^{-1}&\Big[g^{ij}-\frac{2\tau}{2b^2\tau-\rho}b^i\,b^j+\frac{4-\rho\tau}{2b^2\tau-\rho}(l^{i}\,b^{j}+b^{i}\,l^{j})\\[5mm] &-\frac{3\rho b^2\tau^3-\rho^2\tau^2-4b^2\tau^2-2\rho\tau+8}{\tau(2b^2\tau-\rho)}l^i\,l^j\Big] \end{split}\end{equation}
where $b$ is the magnitude of the vector $b^{i}=g^{ij}b_{j}$.

\noindent We obtained the relation between the Cartan connection coefficients $F^{\,i}_{jk}$ and $^{*}F^{\,i}_{jk}$ as \cite{GuPa13}
\begin{equation}^{*}F^{\,i}_{jk}=F^{\,i}_{jk}+D^{\,i}_{jk}\,.\end{equation}
Transvecting by $y^{j}$ and using $F^{\,i}_{jk}\,y^{j}=G^{\,i}_{k}$, we get
\begin{equation}^{*}G^{\,i}_{k}=G^{\,i}_{k}+D^{\,i}_{0k}\,,\end{equation}
where the subscript `$0$' denotes the contraction by the supporting element $y^j$. Further transvecting (3.9) by $y^{k}$ and using $G^{\,i}_{k}\,y^{k}=2\,G^{\,i}$, we get
\begin{equation}2\,^{*}G^{\,i}=2\,G^{\,i}+D^{\,i}_{00}\,.\end{equation}
Differentiating (3.9) partially with respect to $y^{h}$ and using $\dot{\partial_{h}}G^{i}_{k}=G^{\,i}_{kh}$, we have
\begin{equation} ^{*}G^{\,i}_{kh}=G^{\,i}_{kh}+\dot{\partial}_{h}D^{\,i}_{0k}\,,\end{equation}
where $G^{\,i}_{kh}$ are the Berwald connection coefficients. The expressions of $D^{\,i}_{00}$, $D^{\,i}_{0j}$ and $D^{\,i}_{jk}$ are given by
\begin{equation}D^{\,i}_{00}=l^i\,\tau\,(A-E^{}_{00})-\frac{2\tau^2}{2-\rho\tau}A\,m^i+\frac{2L\tau}{2-\rho\tau}\Big(\frac{1}{\beta}\,\beta^{}_{0}\,m^i-F^{\,i}_{0}\Big),\end{equation}
\begin{equation}D^{\,i}_{0j}=l^i\Big\{\frac{1}{\tau}\Big(\frac{2}{\beta}b^2-\frac{\rho}{L}\Big)^{-1}G_{\beta j}+\frac{1}{\tau}G^{}_{j}\Big\}-\frac{2\,m^i}{2-\rho\tau}\Big(\frac{2}{\beta}b^2-\frac{\rho}{L}\Big)^{-1}G_{\beta j}+\frac{L}{2\tau-\rho\tau^2}G^{\,i}_{j}\,, \end{equation}
\begin{equation}D^{j}_{ik}=l^j\Big\{\frac{1}{\tau}\Big(\frac{2}{\beta}b^{2}-\frac{\rho}{L}\Big)^{-1}H^{}_{\beta ik}+\frac{1}{\tau}H_{ik}\Big\}-\frac{2\,m^j}{2-\rho\tau}\Big(\frac{2}{\beta}b^{2}-\frac{\rho}{L}\Big)^{-1}H^{}_{\beta ik}+\frac{L}{2\tau-\rho\tau^2}H^{\,j}_{ik}\,,\end{equation}
where
\begin{equation}\begin{split}G_{ij}=&\frac{\tau^2}{\beta}(m_i\beta_j-m_j\beta_i)-\tau^2\,F_{ij}-\frac{1}{2}(2\tau-\rho\tau^2)L_{ijr}D^{\,r}_{00}-\frac{\tau^2}{2}\rho^{}_{0}L_{ij}\\[3mm]
&-\frac{\tau}{\beta}D^{\,r}_{00}\,\sum_{ijr}m_i\Big\{(\rho\tau-1)L_{jr}-\frac{\tau}{\beta}m_jb_r\Big\}+\frac{\tau}{\beta}\beta^{}_{0}\Big((\rho\tau-1)L_{ij}-\frac{3\tau}{\beta}m_im_j\Big)\,,\end{split}\end{equation}
\begin{equation}G_{j}=\tau^2(F_{j0}-E_{j0})\,,\end{equation}
\begin{equation}\begin{split} H^{}_{jik}&=\frac{(\rho\tau^2-2\tau)}{2}\,\big\{L^{}_{ijr}D^{r}_{0k}+L^{}_{jkr}D^{r}_{0i}-L^{}_{kir}D^r_{0j}\big\}\\
&-\frac{\tau}{\beta}D^{r}_{0k}\,\sum_{ijr}m_i\Big\{(\rho\tau-1)L_{jr}-\frac{\tau}{\beta}m_jb_r\Big\}-\frac{\tau}{\beta}D^{r}_{0i}\,\sum_{jkr}m_j\Big\{(\rho\tau-1)L_{kr}-\frac{\tau}{\beta}m_kb_r\Big\}\\
&+\frac{\tau}{\beta}D^{r}_{0j}\,\sum_{kir}m_k\Big\{(\rho\tau-1)L_{ir}-\frac{\tau}{\beta}m_ib_r\Big\}-\frac{\tau^2}{2}(\rho^{}_{k}L_{ij}+\rho^{}_{i}L_{jk}-\rho^{}_{j}L_{ki})\\
&+\frac{\tau}{\beta}\Big\{\beta^{}_{k}\big((\rho\tau-1)L_{ij}-\frac{3\tau}{\beta}m_im_j\big)+\beta^{}_{i}\big((\rho\tau-1)L_{jk}-\frac{3\tau}{\beta}m_jm_k\big)\\
&-\beta^{}_{j}\big((\rho\tau-1)L_{ki}-\frac{3\tau}{\beta}m_km_i\big)\Big\},
\end{split}\end{equation}
\begin{equation}H^{}_{ik}=\frac{\tau^2}{\beta}(m_i\beta^{}_{k}+m^{}_{k}\beta^{}_{i})-\tau^2E_{ik}-\frac{1}{2}(G_{ik}+G_{ki})\,,\end{equation}
and
\begin{equation}\begin{split}
A=&\Big(\frac{2}{\beta}\beta^{}_{0}m^2-2F^{}_{\beta 0}\Big)\Big(\frac{2}{\beta}b^2-\frac{\rho}{L}\Big)^{-1},\\
F^{\,i}_{j}=&g^{ir}F_{rj}\,,\quad G^{\,i}_j=g^{ik}G_{kj}\,,\quad H^{\,j}_{ik}=g^{jm}H^{}_{mik}\,,\\
2E_{ij}=&b_{i|j}+b_{j|i}\,,\quad 2F_{ij}=b_{i|j}-b_{j|i}\,,\quad \beta^{}_{j}=\beta^{}_{|j}\,,\quad \rho^{}_{k}=\rho^{}_{|k}=\partial_k\rho.\\
\end{split}\end{equation}
\[\textrm{The symbol } \sum^{}_{ijk} \textrm{denote cyclic interchange of indices $i,j,k$ and summation.~~~~~~~~~~~~~~~~~~~~~~~~~~}\]
\begin{lemma}\label{lem4}
\emph{\textbf{\cite{GuPa13}}} If the \textsl{h}-vector $b_i$ is gradient then the scalar $\rho$ is constant.\end{lemma}

\begin{lemma}\label{lem5}
\textbf{\emph{\cite{GuPa13}}} For the Kropina change with \textsl{h}-vector, the difference tensor $D^{\,i}_{jk}$ vanishes if and only if the vector $b_{i}$ is parallel with respect to the Cartan connection of $F^{n}$, and then the Berwald connection coefficients for both the spaces $F^{n}$ and $^{*}F^{n}$ are the same.
\end{lemma}

\section{Hypersurface $^{*}F^{n-1}$ of the space $^{*}F^{n}$}
~~~Let us consider a Finslerian hypersurface $F^{n-1}=(M^{n-1},\underline{L}(u,v))$ of $F^{n}$ and a Finslerian hypersurface $^{*}F^{n-1}=(M^{n-1},\,^{*}\underline{L}(u,v))$ of $^{*}F^{n}$. Let $N^{i}$ be the unit normal vector at a point of $F^{n-1}$. The functions $B^{\,i}_{\alpha}(u)$ may be considered as the components of $n$-1 linearly independent vectors tangent to $F^{n-1}$ and they are invariant under the Randers conformal change. Then the unit normal vector $^{*}N^{i}(u,v)$ of $^{*}F^{n-1}$ is uniquely determined by
\begin{equation} ^{*}g_{ij}\,B^{\,i}_{\alpha}\,^{*}N^{\,j}=0\,,\qquad ^{*}g_{ij}\,^{*}N^{\,i}\,^{*}N^{\,j}=1\,.\end{equation}
The inverse projection factors $^{*}B^{\alpha}_{i}$ are uniquely defined along $^{*}F^{n-1}$ by
\begin{equation} ^{*}B^{\alpha}_{i}=\,^{*}g^{\alpha\beta}\,^{*}g_{ij}\,B^{j}_{\beta}\,,\end{equation}
where $^{*}g^{\alpha\beta}$ is the inverse of the metric tensor $^{*}g_{\alpha\beta}$ of $^{*}F^{n-1}$.\\
From (4.2), it follows that
\begin{equation}B^{i}_{\alpha}\,^{*}B^{\beta}_{i}=\delta^{\beta}_{\alpha}\,,
\quad B^{i}_{\alpha}\,^{*}N_{i}=0\,,\quad ^{*}N^{i}\,^{*}B^{\alpha}_{i}=0\,,\quad ^{*}N^{i}\,^{*}N_{i}=1\,.\end{equation}
Transvecting (2.3) by $v^{\alpha}$, we get
\begin{equation} y^{}_{j}\,N^{j}_{}=0.\end{equation}
Transvecting (3.4) by $N^{i}\,N^{j}$ and paying attention to (2.3) and (4.4), we have
\begin{equation} ^{*}g_{ij}\,N^{i}\,N^{j}=(2\tau^2-\rho\tau^3)+3\tau^4\,(b_{i}\,N^{i})^{2},\end{equation}
which shows that $N^{j}/\sqrt{(2\tau^2-\rho\tau^3)+3\tau^4\,(b_{i}\,N^{i})^{2}}$ is a unit vector. Again transvecting (3.4) by $B^{i}_{\alpha}\,N^{j}$ and using (2.3) and (4.4), we get
\begin{equation} ^{*}g_{ij}\,B^{i}_{\alpha}\,N^{j}=(b_{j}\,N^{j})(3\tau^4\,b_{i}\,B^{i}_{\alpha}-4\tau^3\,l_{i}\,B^{i}_{\alpha}).\end{equation}
This shows that the vector $N^{j}$ is normal to $^{*}F^{n-1}$ if and only if
\[(b_{j}\,N^{j})(3\tau^4\,b_{i}\,B^{i}_{\alpha}-4\tau^3\,l_{i}\,B^{i}_{\alpha})=0\,.\]
This implies at least one of the following :
\[(a)\,\,\, 3\tau^4\,b_{i}\,B^{i}_{\alpha}-4\tau^3\,l_{i}\,B^{i}_{\alpha}=0\,,\qquad\qquad(b)\,\,\,b_{j}\,N^{j}=0\,.\]
The transvection of $3\tau^4\,b_{i}\,B^{i}_{\alpha}-4\tau^3\,l_{i}\,B^{i}_{\alpha}=0$ by $v^{\alpha}$ gives $3\tau^4\,b_{i}\,y^{i}-4\tau^3\,l_{i}\,y^{i}=0$, i.e. $\tau=L/\beta=0$, which is not possible. Therefore $3\tau^4\,b_{i}\,B^{i}_{\alpha}-4\tau^3\,l_{i}\,B^{i}_{\alpha}\neq0$. Hence
\begin{equation} b_{j}\,N^{j}=0\,.\end{equation}
This shows that the vector $N^{j}$ is normal to $^{*}F^{n-1}$ if and only if $b_{j}$ is tangent to $F^{n-1}$.

\noindent From (4.5), (4.6) and (4.7), we can say that $N^{j}/\sqrt{(2\tau^2-\rho\tau^3)}$ is a unit normal vector of $^{*}F^{n-1}$. Therefore in view of (4.1), we get
\begin{equation} ^{*}N^{i}=N^{i}/\sqrt{(2\tau^2-\rho\tau^3)}\,,\end{equation}
which, in view of (3.4), (4.4) and (4.7), gives
\begin{equation} ^{*}N_{i}=\,^{*}g_{ij}\,\,^{*}N^{j}=\sqrt{(2\tau^2-\rho\tau^3)}\,N_{i}\,.\end{equation}
Thus, we have:
\begin{theorem}\label{th3}
Let $^{*}F^{n}$ be the Finsler space obtained from $F^{n}$ by a Kropina change \emph{(1.3)} with \textsl{h}-vector. If $^{*}F^{n-1}$ and $F^{n-1}$ are the hypersurfaces of these spaces then the vector $b_{i}$ is tangential to the hypersurface $F^{n-1}$ if and only if every vector normal to $F^{n-1}$ is also normal to $^{*}F^{n-1}$. \end{theorem}

In view of (4.4) and (4.7), the vector $m_{i}$ appearing in (3.2) satisfies
\begin{equation}m_{i}\,N^{i}=0\,.\end{equation}
As $h_{ij}=g_{ij}-l_{i}l_{j}$, equations (2.3) and (4.4) yield
\begin{equation}h_{ij}\,B^{i}_{\alpha}\,N^{i}=0\,.\end{equation}
Transvecting (3.6) by $B^{i}_{\alpha}\,B^{j}_{\beta}\,N^{k}$ and using (4.10) and (4.11), we get
\begin{equation} ^{*}C_{ijk}\,B^{i}_{\alpha}\,B^{j}_{\beta}\,N^{k}=(2\tau^2-\rho\tau^3)\,C_{ijk}\,B^{i}_{\alpha}\,B^{j}_{\beta}\,N^{k}\,. \end{equation}
Using (2.10) and (4.8), equation (4.12) may be written as
\begin{equation}^{*}M_{\alpha\beta}=\sqrt{(2\tau^2-\rho\tau^3)}\,M_{\alpha\beta}\,.\end{equation}
From (2.8), (3.13) and (4.9), we get
\[^{*}H_{\alpha}=\sqrt{(2\tau^2-\rho\tau^3)}(H_{\alpha}+N_{i}D^{i}_{0k}B^{k}_{\alpha})\,.\]
Transvecting by $v^{\alpha}$ and using $v^{\alpha}\,B^{k}_{\alpha}=y^{k}$, we get
\begin{equation}^{*}H_{0}=\sqrt{(2\tau^2-\rho\tau^3)}(H_{0}+N_{i}D^{i}_{00})\,.\end{equation}
Transvecting (3.12) by $N_{i}$ and using (4.4) and (4.10), we get
\begin{equation}D^{i}_{00}\,N_{i}=-\frac{2L\tau}{2-\rho\tau}\,F^{i}_{0}\,N_{i}\,.\end{equation}
If the vector $b_{i}$ is gradient, i.e. $b_{i|j}=b_{j|i}\,$, then
\begin{equation}F_{ij}=0,\end{equation}
and then by Lemma \ref{lem4}, we have
\begin{equation}\rho_{i}=0.\end{equation}
Therefore (4.15) becomes $D^{i}_{00}\,N_{i}=0\,,$ and then equation (4.14) reduces to
\[^{*}H_{0}=\sqrt{(2\tau^2-\rho\tau^3)}H_{0}\,.\]
Thus, in view of Lemma \ref{lem1}, we have:
\begin{theorem}\label{th4}
Let the \textsl{h}-vector $b_{i}(x,y)$ be a gradient and tangent to the hypersurface $F^{n-1}$. Then the hypersurface $F^{n-1}$ is a hyperplane of the first kind if and only if the hypersurface $^{*}F^{n-1}$ is a hyperplane of the first kind. \end{theorem}

Taking the relative \textsl{h}-covariant differentiation of (4.7) with respect to the Cartan connection of $F^{n-1}$, we get
\[b_{i|\beta}\,N^{i}+b_{i}\,N^{i}_{\,\,|\beta}=0.\]
Using (2.11) and (2.12), the above equation gives
\[(b_{\,i|j}\,B^{j}_{\beta}+b_{i}|_{j}\,N^{j}\,H_{\beta})N^{i}-b_{i}\,H_{\alpha\beta}\,B^{\alpha}_{j}\,g^{ij}=0.\]
Transvecting by $v^{\beta}$ and using (2.9), we get
\begin{equation}b_{i|0}\,N^{i}=(H_{\alpha}+M_{\alpha}\,H_{0})B^{\alpha}_{j}\,b^{j}-b_{i}|_{j}\,H_{0}\,N^{i}\,N^{j}. \end{equation}
For a hypersurface of the first kind, $H_0=0=H_\alpha\,$. Then (4.18) reduces to $b_{i|0}\,N^{i}=0$. If the vector $b_{i}$ is gradiant, i.e. $b_{i|j}=b_{j|i}\,$, then we get
\begin{equation}E_{i0}\,N^{i}=b_{i|0}\,N^{i}=\beta^{}_{i}\,N^i=0.\end{equation}
Therefore (3.16) gives
\begin{equation}G_{j}N^{j}_{}=0\,.\end{equation}
Transvecting (3.15) by $b^iN^j$ and using (4.10), (4.16), (4.17) and (4.19), we get
\begin{equation}G_{ij}\,b^i\,N^{j}_{}=G^{}_{\beta j}\,N^{j}_{}=0\,.\end{equation}
Again transvecting (3.15) by $N^iB^{j}_{\alpha}$ and using $D^{\,i}_{00}N^{}_{i}=0$ and $b^{}_{i}N^i=0$, we obtain
\begin{equation}G_{ij}\,N^iB^{j}_{\alpha}=0\,.\end{equation}
Transvecting (3.13) by $N^{}_{i}B^{j}_{\alpha}$ and using (4.4), (4.10) and (4.22), we get
\begin{equation}D^{\,i}_{0j}\,N^{}_{i}B^{j}_{\alpha}=0\,.\end{equation}
Transvecting (3.12) by $L_{ijk}$ and using (3.5) and (4.16), we have
\begin{equation}L_{ijk}D^{\,i}_{00}=(A-\beta^{}_{0})\Big[-\Big(\frac{\tau}{L^2}+\frac{4\,\rho\tau^2}{L^2(2-\rho\tau)}\Big)\,h^{}_{jk}+\frac{2\tau^2}{L^2(2-\rho\tau)}(m_il_k+m_kl_i)\Big]\,.\end{equation}
Transvecting (3.13) by $N^{j}_{}B^{k}_{\alpha}\,h^{}_{ki}$ and using (4.10) and (4.24), we get
\begin{equation}D^{\,i}_{0j}N^{j}_{}B^{k}_{\alpha}\,h^{}_{ik}=0\,.\end{equation}
Again transvecting (3.13) by $b^{}_{i}N^{j}_{}$ and using (4.20) and (4.21), we get
\begin{equation}D^{\,i}_{0j}\,b^{}_{i}N^{j}_{}=0\,.\end{equation}
Transvecting (3.17) by $N^jB^{i}_{\alpha}B^{k}_{\beta}$ and using (4.10), (4.23), (4.25) and (4.26), we obtain
\begin{equation}H^{}_{jik}N^jB^{i}_{\alpha}B^{k}_{\beta}=\frac{(\rho\tau^2-2\tau)}{2}N^jB^{i}_{\alpha}B^{k}_{\beta}\big(L^{}_{ijr}D^{\,r}_{0k}+L^{}_{jkr}D^{\,r}_{0i}-L^{}_{kir}D^{\,r}_{0j}\big)\,.\end{equation}
Let the \textsl{h}-vector $b_i$ satisfies the condition
\begin{equation}b^{}_{r|0}\,C^{\,r}_{ij}=0\,,\end{equation}
then
\begin{equation}\beta^{}_{r}\,C^{\,r}_{ij}=0\,.\end{equation}
Transvecting $G^{}_{rk}$ by $C^{\,r}_{ij}\,N^jB^{i}_{\alpha}B^{k}_{\beta}$ and using (4.16), (4.17), (4.24) and (4.29), we get
\[\begin{split}C^{\,r}_{ij}\,N^jB^{i}_{\alpha}B^{k}_{\beta}\,G^{}_{rk}=-C^{\,r}_{ij}\,N^jB^{i}_{\alpha}&B^{k}_{\beta}\Big[(A-\beta^{}_{0})\Big(\frac{\tau}{L^2}+\frac{4\,\rho\tau^2}{L^2(2-\rho\tau)}\Big)\,h^{}_{rk}\\
&+\frac{\tau}{\beta}(\rho\tau-1)D^{\,s}_{00}m^{}_{s}L^{}_{rk}-\frac{\tau}{\beta}\beta^{}_{0}(\rho\tau-1)L^{}_{rk}\Big]\,.\end{split}\]
Using (2.10), above equation is simplified as
\[\begin{split}C^{\,r}_{ij}\,N^jB^{i}_{\alpha}B^{k}_{\beta}\,G^{}_{rk}=\lambda\,M^{}_{\alpha\beta}\,,\end{split}\]
where we put
\[\lambda=-(A-\beta^{}_{0})\Big(\frac{\tau}{L^2}+\frac{4\,\rho\tau^2}{L^2(2-\rho\tau)}\Big)+(A-\beta^{}_{0})(\rho\tau-1)\,\frac{2\tau^3m^2}{L\beta\,(2-\rho\tau)}+\frac{\tau}{L\beta}\,\beta^{}_{0}(\rho\tau-1).\]
Hence we get
\begin{equation}L^{}_{ijr}N^{j}_{}B^{i}_{\alpha}B^{k}_{\beta}\,D^{\,r}_{0k}=\frac{2\,\lambda}{2\tau-\rho\tau^2}\,M^{}_{\alpha\beta}\,.\end{equation}
By interchanging the indices $i$ and $k$, we have
\begin{equation}L^{}_{jkr}N^{j}_{}B^{i}_{\alpha}B^{k}_{\beta}\,D^{\,r}_{0i}=\frac{2\,\lambda}{2\tau-\rho\tau^2}\,M^{}_{\alpha\beta}\,.\end{equation}
Transvecting (3.15) by $N^{j}_{}$ and using (4.10), (4.16), (4.19) and (4.24), we get
\[G^{}_{ij}N^{j}_{}=\mu\,N^{}_{i}\,,\]
where we put
\[\mu=\frac{1}{2}(2\tau-\rho\tau^2)(A-\beta^{}_{0})\Big(\frac{\tau}{L^2}+\frac{4\,\rho\tau^2}{L^2(2-\rho\tau)}\Big)-\frac{\tau}{L\beta}m^{}_{s}D^{\,s}_{00}+\frac{\tau}{L\beta}\beta^{}_{0}(\rho\tau-1).\]
Using this, we obtain
\begin{equation}L^{}_{kir}N^{j}_{}B^{i}_{\alpha}B^{k}_{\beta}\,D^{\,r}_{0j}=\frac{2\mu}{2\tau-\rho\tau^2}\,M^{}_{\alpha\beta}\,.\end{equation}
Using (4.30), (4.31), (4.32) in (4.27), we get
\begin{equation}H^{}_{jik}N^jB^{i}_{\alpha}B^{k}_{\beta}=(\mu-2\lambda)M^{}_{\alpha\beta}\,.\end{equation}
Since $l^{j}_{}N_j=0$ and $m^{j}_{}N_j=0$, the equation (3.14) gives
\begin{equation}D^{\,j}_{ik}N^{}_{j}B^{i}_{\alpha}B^{k}_{\beta}=\frac{(\mu-2\lambda)L}{2\tau-\rho\tau^2}M^{}_{\alpha\beta}\,.\end{equation}
From (2.7), (3.8) and (4.9), we have
\begin{equation}^{*}H_{\alpha\beta}-\,^{*}M_{\alpha}\,^{*}H_{\beta}=\sqrt{2\tau^2-\rho\tau^3}\,(H_{\alpha\beta}+D^{i}_{jk}\,N_{i}\,B^{j}_{\alpha}\,B^{k}_{\beta})-M_{\alpha}\,H_{\beta}\,.\end{equation}

\noindent Thus from (4.28), (4.34) and (4.35), we have:
\begin{theorem}\label{th5}
For the Kropina change with \textsl{h}-vector, let the \textsl{h}-vector $b_{i}$ be a gradient and tangential to the hypersurface $F^{n-1}$ and satisfies the condition \emph{(4.28)}. Then
\begin{enumerate}
  \item $^{*}F^{n-1}$ is a hyperplane of the second kind if $F^{n-1}$ is a hyperplane of the second kind and $M_{\alpha\beta}=0$.
  \item $^{*}F^{n-1}$ is a hyperplane of the third kind if $F^{n-1}$ is a hyperplane of the third kind.
\end{enumerate}\end{theorem}

A Finsler space $F^{n}$ is called a Landsberg space if the (\textsl{v})\textsl{hv}-torsion tensor $P^{\,r}_{ij}$ vanishes, i.e.
\begin{equation}P^{\,r}_{ij}=C^{\,r}_{ij\,|k}\,y^k=0.\end{equation}
Taking h-covariant derivative of (1.3)$(ii)$ and using $L^{}_{|k}=0=h^{}_{ij|k}$, (4.17), we get
\[b^{}_{r|k}\,C^{\,r}_{ij}+b^{}_{r}\,C^{\,r}_{ij|k}=0.\]
Contracting by $y^k$ and using (4.36), we get
\[b^{}_{r|0}\,C^{\,r}_{ij}=0,\]
which is the condition (4.28). Hence we have:
\begin{theorem}\label{th6}
For the Kropina change with \textsl{h}-vector, let the \textsl{h}-vector $b_{i}$ be a gradient and tangential to the hypersurface $F^{n-1}$ of a Landsberg space $F^{n}$. Then
\begin{enumerate}
  \item $^{*}F^{n-1}$ is a hyperplane of the second kind if $F^{n-1}$ is a hyperplane of the second kind and $M_{\alpha\beta}=0$.
  \item $^{*}F^{n-1}$ is a hyperplane of the third kind if $F^{n-1}$ is a hyperplane of the third kind.
\end{enumerate}\end{theorem}

For the Kropina change with \textsl{h}-vector, let the vector $b_{i}$ be parallel with respect to the Cartan connection of $F^{n}$. In view of Lemma \ref{lem5}, we have
\begin{equation}^{*}F^{\,i}_{jk}=F^{\,i}_{jk}\,.\end{equation}
Thus from (2.7), (4.9), (4.13) and (4.37), we have:
\begin{theorem}\label{th7}
For the Kropina change with \textsl{h}-vector, let the vector $b_{i}$ be parallel with respect to the Cartan connection of $F^{n}$ and tangent to the hypersurface $F^{n-1}$. Then $^{*}F^{n-1}$ is a hyperplane of the second \emph{(}third\emph{)} kind if and only if $F^{n-1}$ is also a hyperplane of the second \emph{(}third\emph{)} kind. \end{theorem}

\end{document}